\documentclass[a4paper,10pt]{article}
\include{amsmath}
\usepackage{amsmath}
\usepackage{amsmath}
\usepackage{amsfonts}
\usepackage{amssymb}

\title{The simplest and fastest method of solving Exact Differential Equation by developing the idea of Basic Functions}
\author {Seyed Ahmad Sabok-Sayr \footnote{e-mail: ahmads@sfu.ca / saboksayr@gmail.com}\\
\normalsize{Departments of Mathematics and Physics}\\
\normalsize{Simon Fraser University}\\
\normalsize{ Burnaby, British Columbia, Canada, V5A 1S6}}
\date{}

\begin{document}

\maketitle
\begin{abstract}
This paper concerns exact differential equations. First, I define two types of functions which I have named Basic Function of Type One and Basic Function of Type Two. I then derive the property and theorems of these functions. Finally, by applying the property and theorems of Basic Functions, I establish the method of solving exact differential equation with n variables which is significantly simpler and faster than the standard method.
\end{abstract}
\section*{Introduction}
From ordinary differential equation, we know the standard method to solve exact differential equations. The hard part of the standard method is the calculation of integration constant and the difficulty increases rapidly as the number of variables increases. In general, for exact differential equation with $n$ variables, we need to find the integration constant of $n-1$ variables which requires tedious calculations. For example, for exact differential equations with $2, 3, 4$ and $5$ variables, the number of integrations and differentiations required to solve them are $5, 17, 49$ and $129$ respectively. The method presented in this paper does not require any calculation of integration constant and to solve the exact differential equation with $n$ variables, in general, we just need $n$ integrations. Therefore, for example for exact differential equation with $5$ variables, using this method, we just need $5$ integrations instead of $129$ integrations and differentiations required in standard method. 
We can get the idea of new method from the standard one. Let's consider the following first-order differential equation :
\begin{equation}
M(x,y) \,\,dx + N(x,y)\,\, dy = 0  \label{EDE2v1}
\end{equation} 
From any ordinary differential equation text book such as \cite[Theorem 2.6.1, p. 91]{boyce}, we know that equation \ref{EDE2v1} is an exact differential equation if and only if the following relation is true:\\
\begin{equation}
M_y (x,y) = N_x(x,y)  \label{Exactness}
\end{equation}
That is , there exists a function $u(x,y)$ such that:
\begin{equation}
 u_x(x,y) = M(x,y),\,\,\,\,\,\,\,\,\,\,\,\,\,\,u_y(x,y) = N(x,y)  \label{ExistFun} 
\end{equation}
From equation \ref{ExistFun} we get:
\begin{equation}
 u(x,y) = \int  M(x,y) dx + h(y)   \label{Soln1}
\end{equation}
Alternatively, we can find $u(x,y)$ as below:
\begin{equation}
 u(x,y) = \int  N(x,y) dy + k(x)   \label{Soln2}
\end{equation}
The idea is that instead of using standard method; using equation \ref{Soln1} or \ref{Soln2}, suppose that the solution $u(x,y)$ is a linear combination of $m$ functions as follows:
\begin{equation}
 u(x,y) = \phi_1 + \phi_2 + ...+ \phi_m  \label{LinearComb}
\end{equation}
where, every function $\phi_i$ can be a function of  $x, y $ or both and $m$ is an integer greater than zero. From equation \ref{Soln1} we conclude that, some of $\phi_is$ are in $\int M(x,y) dx$ and the rest of $\phi_is$ are in $h(y)$. Alternatively, from equation \ref{Soln2}, we conclude that some of the $\phi_is$ are in $\int N(x,y) dy$ and the rest of $\phi_is$ are in $k(x)$. Instead of using equations \ref{Soln1} or \ref{Soln2} separately, let's consider the union of the $\phi_is$ which are in $\int M(x,y) dx$ and $\phi_is$ in $\int N(x,y) dy$. We may ask the following questions: Does the union contain all $\phi_is$ ? Are any $\phi_is$ missed in the union? In this paper we show  how to derive the solution out of the union of $\phi_is$. We also generalize the method for exact differential equation with $n$ variables.

\subsection*{Theory}
\textsc{Definition 1.} Let $\Omega$ be a nonempty simply-connected open region of $\mathbb{R}^n$, and $X=(x_1,x_2,...,x_n) \in \Omega$. Suppose that $f :\Omega\rightarrow\mathbb{R}$ be a function such that $\frac{\partial f}{\partial x_i}$ exists and is integrable in $\Omega$ for some $i$ such that $1 \leq i \leq n$. The function $f$ is called a \textbf{Basic Function Of Type One} with respect to $x_i$ if and only if the following condition is satisfied when we set the integration constant to zero:
\begin{equation} 
\int\dfrac{\partial}{\partial x_i} f(x_1,\ldots,x_i,\ldots,x_n) dx_i = f(x_1,\ldots,x_i,\ldots,x_n)  \label{Defn1}
\end{equation}

\textsc{Definition 2.} If a function is Basic Function of Type One with respect to all of its variables, then it is a \textbf{Basic Function of Type Two}.\\

\textsc{Note 1.} If a Basic Function of Type One has just one variable, by definition $2$, it is Basic Function of Type Two as well.\\

\textsc{Example 1}: Let $f(x,y) = \sin(x)\,\cos(y)+\,\sin(y)$ then:\\
\begin{eqnarray*}
   \int\dfrac{\partial}{\partial x} f(x,y) dx &=& \sin(x)\,\cos(y) + \overbrace{C_1(y)}^{set\,\,to\,\,zero} \\
&=& \sin(x)\,\cos(y)\\
&\neq& f(x,y) \\
\end{eqnarray*}
\begin{eqnarray*}
 \int\dfrac{\partial}{\partial y} f(x,y) dy &=& \sin(x)\,\cos(y) +\,\sin(y)+ \overbrace{C_2(x)}^{set\,\,to\,\,zero} \\
&=& \sin(x)\,\cos(y) +\,\sin(y)\\
&=&f(x,y) 
\end{eqnarray*}
\\
so $f(x,y)$ is not a Basic Function of Type One with respect to $x$ but it is a Basic Function of Type One with respect to $y$  .\\

\textsc{Example 2.}  Let $g(x,y) = \sin(x)\, \cos(y)$ , then:\\
\begin{eqnarray*}
  \int\dfrac{\partial}{\partial x} g(x,y) dx &=&\sin(x)\, \cos(y) + \overbrace{C_3(y)}^{set\,\,to\,\,zero}\\
&=& \sin(x)\, \cos(y)\\
&=& g(x,y)
\end{eqnarray*}
\begin{eqnarray*}
\int\dfrac{\partial}{\partial y} g(x,y) dy &=& \sin(x)\, \cos(y) + \overbrace{C_4(x)}^{set\,\,to\,\,zero} \\
&=&\sin(x)\, \cos(y)\\
&=& g(x,y)
\end{eqnarray*}
\\
hence $g(x,y)$ is a Basic Function of Type One with respect to both $x$ and $y$, so it is a Basic Function of Type Two . \\

\textsc{Property 1.} Suppose that $B(x_1,\ldots,x_n)$ is a basic function of type two. let's consider that we are given $dB$ as below:
\begin{eqnarray*}
dB = \frac{\partial B}{\partial x_1} dx_1+\ldots+\frac{\partial B}{\partial x_n} dx_n
\end{eqnarray*}
We can derive $B$ from $dB$ by integrating from $dB$ and eliminating the repeated terms as below:
\begin{eqnarray*}
\int\frac{\partial B}{\partial x_1} dx_1+\ldots+\int\frac{\partial B}{\partial x_n} dx_n=B+ \overbrace{B +\ldots+ B}^{to\,be\,eliminated}= B(x_1,\ldots,x_n)\\
\end{eqnarray*}
In fact, when we integrate from $dB$, we produce as many $B$s as the number of variables of $B$.\\

\textsc{Example 3.} Let's consider the function $g(x,y)=\sin(x)\cos(y)$ which we used in example $2$. Now, we know that it is a Basic Function of Type Two. Suppose that the differential of $g(x,y)$ is given as below:\\
\begin{eqnarray*}
 dg(x,y) &=& \dfrac{\partial g}{\partial x} dx + \dfrac{\partial g}{\partial y} dy \\
&=&\cos(x)\cos(y) dx - \sin(x)\sin(y) dy\\
\end{eqnarray*}
To solve the above equation, we only need to apply the property $1$ ; integrate each part, and eliminate repeated terms as below:\\
\begin{eqnarray*}
 g(x,y) &=& \int \cos(x)\cos(y) dx - \int \sin(x)\sin(y) dy \\
&=& \sin(x)\cos(x) + \overbrace{\sin(x) \cos(x)}^{to\, be\, eliminated} \\
&=& \sin(x) \cos(x) 
\end{eqnarray*}

\textsc{Theorem 1.} Let $\Omega$ be a nonempty simply-connected open region of $\mathbb{R}^n$, and $X=(x_1,x_2,...,x_n) \in \Omega$. Suppose that $f :\Omega\rightarrow\mathbb{R}$ be a function such that $\frac{\partial f}{\partial x_i}$ exists and is integrable in $\Omega$ for some $i$ such that $1 \leq i \leq n$. Then $\int \frac{\partial}{\partial x_i} F(x_1,\ldots,x_n) dx_i$ is a Basic Function of Type One with respect to $x_i$.\\

\textit{Proof.} If $\int \frac{\partial}{\partial x_i} F(x_1,\ldots,x_n) dx_i$ is a Basic Function of Type One, it must satisfy definition $1$. Let  $g(x_1,\ldots,x_n)=\int \frac{\partial}{\partial x_i} F(x_1,\ldots,x_n) dx_i$, and apply definition $1$ for $g$:
\begin{eqnarray*}
\int \dfrac{\partial}{\partial x_i} g(x_1,\ldots,x_n) dx_i & = & \int \left[\dfrac{\partial}{\partial x_i} \int \dfrac{\partial } {\partial x_i} F(x_1,\ldots,x_n)dx_i\right]dx_i \\
& = & \int \dfrac{\partial}{\partial x_i} F(x_1,\ldots,x_n) dx_i \\
& = & g(x_1,\ldots,x_n)
\end{eqnarray*}
hence $\int \frac{\partial}{\partial x_i} F(x_1,\ldots,x_n) dx_i$ is a Basic Function of Type One with respect to $x_i$.\\

\textsc{Example 4.} Let's consider the function $f(x,y)=\sin(x)\cos(y) + \sin(y)$ which we used in example $1$. Now, we know that it is not a Basic Function of Type One with respect to $x$, but $\int \frac{\partial f}{\partial x} dx = \sin(x) \cos(y)$ which is a Basic Function of Type One with respect to $x$.\\
\begin{eqnarray*}
\int \frac{\partial f}{\partial x} dx = \sin(x) \cos(y)
\end{eqnarray*}
\\
Using theorem $1$, we can make a Basic Function of Type One from a given function.\\

\textsc{Theorem 2.} Let $\Omega$ be a nonempty simply-connected open region of $\mathbb{R}^n$, and $X = (x_1,x_2,...,x_n) \in \Omega$. Suppose that $f :\Omega\rightarrow\mathbb{R}$ be a function such that $\frac{\partial f}{\partial x_i}$ exists and is integrable in $\Omega$ for some $i$ such that $1 \leq i \leq n$, then $F(x_1,\ldots,x_n)$ is a linear combination of some Basic Functions of Type Two and some constant.\\

\textit{Proof:} Let's check the definition of Basic Function from $x_1$ to $x_n$:
\begin{eqnarray*}
  \int \frac{\partial }{\partial x_1} F(x_1,\ldots,x_n) dx_1 = f_1 (x_1,\ldots, x_n) 
\end{eqnarray*}
By theorem 1, $f_1$ is a Basic Function of Type One with respect to $x_1$. Now differentiate with respect to $x_1$ from both sides,  then we get:
\begin{eqnarray*}
 \frac{\partial}{\partial x_1} F(x_1,\ldots,x_n) = \frac{\partial}{\partial x_1}f_1 (x_1,\ldots,x_n) \label{Thm2}
\end{eqnarray*}
It follows that:
\begin{eqnarray*}
F(x_1,\ldots,x_n) - f_1(x_1,\ldots,x_n) = C_1(x_2,\ldots,x_n)
\end{eqnarray*}
where $C_1(x_2,\ldots,x_n)$ is a constant with respect to $x_1$ but in general it can be a function of $x_2,\ldots,x_n$. Now let's work on $C_1$ as below:
\begin{eqnarray*}
 \int \frac{\partial}{\partial x_2}C_1(x_2,\ldots,x_n) dx_2 = f_2(x_2,\ldots,x_n)
\end{eqnarray*}
By theorem 1, $f_2$ is a Basic Function of Type One with respect to $x_2$. Let's differentiate with respect to $x_2$ from both sides, then we get:
\begin{eqnarray*}
 \frac{\partial}{\partial x_2} C_1(x_2,\ldots,x_n) = \frac{\partial}{\partial x_2}f_2 (x_2,\ldots,x_n)
\end{eqnarray*}
It follows that:
\begin{eqnarray*}
C_1(x_2,\ldots,x_n) - f_2(x_2,\ldots,x_n) = C_2(x_3,\ldots,x_n)
\end{eqnarray*}
Therefore:
\begin{eqnarray*}
 F(x_1,\ldots,x_n) - f_1 (x_1,\ldots,x_n) - f_2(x_2,\ldots,x_n) = C_2(x_3,\ldots,x_n)
\end{eqnarray*}
If we keep doing the above process for $C_2,\ldots, C_n$ we get the following result.
\begin{eqnarray*}
 F(x_1,\ldots,x_n) - f_1(x_1,\ldots,x_n) - \ldots - f_n(x_n) = C_n
\end{eqnarray*}
It follows that:
\begin{equation}
 F(x_1,\ldots,x_n) = \sum _{i=1}^n f_i + C
\end{equation}
where $f_i$ is a Basic Function of Type One with respect to $x_i$ for $i=1,\ldots,n$ and $C = C_n$ is a constant. We can apply this result; equation \ref{Thm2} to all $f_i$ and keep doing that for the following Basic Functions of Type One. Since we can not do that infinitely many times, this process will stop at some point. Therefore, the final Basic Functions all have either one variable or they are Basic Functions of Type Two. So all final $f_i$s are Basic Functions of Type Two.\\

\textsc{Example 5.} Again, let's consider the function $f(x,y)=\sin(x)\cos(y) + \sin(y)$ which we used in example $1$. Now, we know that it is a Basic Function of Type One with respect to $y$, but it is not a Basic Function of Type One with respect to $x$. It is obvious that the $f(x,y)$ is sum of $f_1 = \sin(x)\cos(y)$ which is a Basic Function of Type Two and $f_2 = \sin(y)$ which is a Basic Function of Type One with respect to $y$ and since it has just one variable it is a Basic Function of Type Two as well. Therefore, $f$ is a sum of two Basic Functions of Type Two.\\
\\
\subsection*{Basic Function Method}
In this section, we apply theorem $2$ and property $1$ to introduce the Basic Function Method. First, Let's consider the exact differential equation with two variables as follows:
\begin{equation}
 M(x,y) dx + N(x,y) dy = 0 \label{EDE2v2}
\end{equation}
Suppose that the solution of the differential equation is $\phi (x,y) = C$  where $C$ is a constant. By theorem $2$,  $\,\phi(x,y)$ is a linear combination of some Basic Functions of Type Two and some constant. Without loss of generality we can ignore the constant. Therefore, the most general form that we can write $\phi(x,y)$ as linear combination of Basic Functions of Type Two is as below:
\begin{equation}
 \phi(x,y) = B_1(x) + B_2(y) + B_3(x,y)  \label{Basic2v}
\end{equation}
where $B_i$ are basic functions of type two. Therefore we can write equation \ref{EDE2v2} as follows:
\begin{center}
$ \left [\dfrac{\partial B_1}{\partial x} + \dfrac{\partial B_3}{\partial x}\right ] dx + \left [\dfrac{\partial B_2}{\partial y} + \dfrac{\partial B_3}{\partial y}\right ] dy = 0$
\end{center}
Then by property $1$ we get:
\begin{eqnarray*}
 \int \left [\dfrac{\partial B_1}{\partial x}  +  \dfrac{\partial B_3}{\partial x}\right ] dx + \int\left [\dfrac{\partial B_2}{\partial y} + \dfrac{\partial B_3}{\partial y}\right ] dy \\
= B_1(x) + B_3(x,y) + B_2(y) + \underbrace{B_3(x,y)}_{repeated}\\
= B_1(x)+B_2(y)+B_3(x,y) = \phi(x,y)\\
\end{eqnarray*}
Which is the solution of exact differential equation \ref{EDE2v2}.\\
\\
\\
Now let generalize the method and consider the exact differential equation with $n$ variables as below:\\
\begin{equation}
 M_1(x_1,\ldots,x_n)dx_1 +M_2(x_2,\ldots,x_n)dx_2 +\ldots+ M_n(x_1,\ldots,x_n)dx_n = 0 \label{EDEnv}
\end{equation}
Suppose the solution of the differential equation is $ \phi (x_1,\ldots,x_n) = C$ where $C$ is a constant. By theorem $2$, $\phi$ is a linear combination of some Basic Functions of Type Two as below:  \\
\begin{equation}
\phi (x_1,\ldots,x_n) = B_1+B_2+\ldots+B_m  \label{Basicnv}
\end{equation}
Where $B_i$  ; $j=1,\ldots,m$ ; are basic Functions of Type Two with respect to some $x_i$ s. Equations \ref{EDEnv} and \ref{Basicnv} implies the following :
\begin{equation}
 d\phi = dB_1+dB_2+\ldots+dB_m    \label{dBasicnv}
\end{equation}
\begin{equation*}
      \,\,\,\,\,\,\,\,\,\,\,\,\,\,\,\,\,\,\,\,\,\,\,\,\,\,\,\, = M_1dx_1 +M_2dx_2 +\ldots+ M_ndx_n   \label{dfunnv}
\end{equation*}
Now, by applying property $1$ to all $B_i$s in equation \ref{dBasicnv}, the solution of equation \ref{EDEnv}, $\phi(x_1,\ldots,x_n)$, is the direct result of the following two rules:\\

\textsc{Rule 1:} Ignoring integration constants, integrate each term of equation \ref{EDEnv} separately as follows :\\
\begin{equation}
\phi = \int M_1dx_1 +\int M_1dx_2 +\ldots+ \int M_ndx_n = C\\   \label{Rule1}
\end{equation}

\textsc{Rule 2:} eliminate any repeated terms in equation \ref{Rule1} after finding the integrals.\\
\\
Therefore, the method of Basic Function can be summarized by Rules $1$ and $2$ which can be applied on an exact differential equation with any number of variables.
\subsection*{An Example Of Exact Differential Equation}

Let's solve the following exact differential equation that motivates the method of Basic Function: 
\begin{eqnarray*}
 (e^x\,\sin \,y\,\cos\, z\, -\, 2y\, \sin\, x\, e^z) dx\, +\,(e^x\,\cos\, y\, \cos\, z\, +\, 2\cos\, x\, e^z\, +\,\frac{1}{yz})dy \,\\
+\, (2 y \,\cos\, x\, e^z\, -\, e^x\, \sin\, y\, \sin\, z \,-\, \frac{\ln \,y}{z^2})dz\, =\, 0
\end{eqnarray*}
where: 
\begin{eqnarray*}
M(x,y,z)&=&e^x\,\sin \,y\,\cos\, z\, -\, 2y\, \sin\, x\, e^z,\\
N(x,y,z)&=&e^x\,\cos\, y\, \cos\, z\, +\, 2\cos\, x\, e^z\, +\,\frac{1}{yz}\\  
P(x,y,z)&=&2 y \,\cos\, x\, e^z\, -\, e^x\, \sin\, y\, \sin\, z \,-\, \frac{\ln \,y}{z^2}
\end{eqnarray*}
The differential equation is exact because $M_y=N_x\;,\;M_z=P_x\;,\;$and$\;N_z=P_y.$
\begin{enumerate}

 \item \textsc{Standard Method.}
\begin{eqnarray*}
 \Phi(x,y,z) & = & \int M(x,y,z) dx = e^x\,\sin\,y\,\cos\,z\,+\,2y\,\cos\,x\,e^z \,+\,h(y,z)\\
\frac{\partial \Phi}{\partial y}\, & = &\,N
\end{eqnarray*}
Thus:
\begin{eqnarray*}
 \,e^x\,\cos\,y\,\cos\,z\,+\,2\cos\,x\,e^z +  \frac{\partial h}{\partial y} = e^x \cos\,y\,\cos\,z + 2\cos\,x\,e^z + \frac{1}{yz}
\end{eqnarray*}
Therefore:
\begin{eqnarray*}
h(y,z)& = & \int \frac{1}{yz} dy = \dfrac{\ln y}{z} + k(z)
\end{eqnarray*}
It follows that:
\begin{eqnarray*}
\Phi (x,y,z) & = & e^x\,\sin\,y\,\cos\,z\,+\,2y\,\cos\,x\,e^z \,+\,\dfrac{\ln y}{z} \,+ \,k(z)
\end{eqnarray*}
We can calculate $k(z)$ as bellow:
\begin{eqnarray*}
\frac{\partial \Phi}{\partial z}\, & = &\,p\,\Rightarrow \, -e^x\sin\,y\,\sin\,z+ 2y\,\cos\,x\,e^z+\ln \, y \,\ln\,z\,+k'(z)\\
& = & -e^x\,\sin\,y\,\sin\,z+2y\,\cos\,x\,e^z + \ln\,y\,\ln\,z
\end{eqnarray*}
Hence:
\begin{eqnarray*}
k'(z) \, & = & \,0 \,\Rightarrow\, k(z) = constant \\
\end{eqnarray*}
we conclude that: 
\begin{eqnarray*}
\Phi(x,y,z) & = &e^x\,\sin\,y\,\cos\,z\,+\,2y\,\cos\,x\,e^z \,+\,\dfrac{\ln y}{z} \, =C\\
\end{eqnarray*}

 \item \textsc{Basic Function Method.} Let's use Rule $1$ and $2$; calculate $\int M dx + \int N \,dy + \int P \,dz$ and eliminate the repeated terms: 
\begin{eqnarray*}
\Phi &=&\,e^x \sin\,y\,\cos\,z\,+\,2y\,\cos\,x\,e^z\,+\underbrace{e^x \sin\,y\,\cos\,z}_{\text {repeated}}\\
&+& \underbrace{2y\,\cos\,x\,e^z}_{\text{repeated}}\,+\frac{\ln\,y}{z}+\,\underbrace{e^x\,\sin\,y\,\cos\,z}_{repeated}\\
&+&\underbrace{2y\,\cos\,x\,e^z}_{repeated} \,+\,\underbrace{\frac{\ln\,y}{z}}_{repeated}\\
 &=& \,e^x \sin\,y\,\cos\,z\,+\,2y\,\cos\,x\,e^z\,+\frac{\ln\,y}{z} = C
\end{eqnarray*}
which gives the same result as the standard method.
\end{enumerate}

 



\end{document}